# THE COEFFICIENTS OF THE REDUCED BARTHOLDI ZETA FUNCTION

MAEDEH S. TAHAEI, SEYED NASER HASHEMI

ABSTRACT. In this paper, we establish a new zeta function based on the Bartholdi zeta function for an undirected graph $G$ called the reduced Bartholdi zeta function. We study the relation between its coefficients and the structure of the graph, and demonstrate that the coefficients count the star subgraphs in the symmetric digraph $\mathcal{D}(G)$. Moreover, we investigate the properties of semi principle minors extracted from the adjacency matrix of the oriented line graph of $G$. We also present a general formula for calculating all the coefficients of the reduced Bartholdi zeta function.

## 1. INTRODUCTION

The Ihara zeta function was first introduced by Ihara [6] in 1966 for finite $k$-regular graphs, in the context of discrete groups of the p-adic zeta functions and it was shown that the Ihara zeta function can be represented in the form of the reciprocal of a polynomial. In 1989, Hashimoto [4] deduced multi variable zeta functions for bi-regular bipartite graphs. Also, a generalization of the determinant expression of the Ihara zeta function to all finite irregular graphs by its adjacency matrix was performed by Bass [2]. Further, Horton [5] demonstrated the relation between the girth of graphs and the Ihara zeta function. Stark and Terras [14, 17] argued about various types of the zeta functions for any graph and presented a comprehensive overview of the zeta functions. In 2008, an investigation about the coefficients of the Ihara zeta functions and their associated cyclical structure in graph has been followed by Scott and Storm [13].

The two-variable zeta function of a graph was first demonstrated by Bartholdi [1]. He obtained a determinant expression for the two-variable zeta function, called the Bartholdi zeta function and gave a formula for concluding the number of bumps on the paths in a graph. Mizuno and Sato [8, 9] defined a new zeta function for a digraph by using a determinant expression of the zeta function and also expressed its weighted Bartholdi zeta function. For generalizing the Bartholdi zeta function, Sato [11], represented a three-variable zeta function of a bipartite graph and presented its determinant expression. Moreover, in [12], he presented an $(n + 1)$-variable Bartholdi L-function of graph $G$. By the weighted scattering matrix of $G$, Oren [10] stated a new version of the Bartholdi zeta function.

There is a substantial research literature on studying the properties of the Ihara zeta function [4, 5, 7], in distinguishing co-spectral graphs [15] and on developing the various types of the Ihara zeta function [3, 16].

---

However, what is still lacking is an explicit description of the relation between the Bartholdi zeta function and the structure of a graph. Furthermore, despite the Ihara zeta function can be computed in polynomial time, the Bartholdi zeta function is not easy to compute, since it is expressible as a two-variable determinant. Note that all types of zeta functions are defined on md2 graphs in which each vertex has at least degree 2.

In this paper, we define a new Bartholdi zeta function of a graph called the reduced Bartholdi zeta function. We investigate about the coefficients of the reduced Bartholdi zeta function and demonstrate that these coefficients can count the special complete bipartite structures as stars in the graph.

The reduced Bartholdi zeta function focuses on the variable related to the number of bumps on paths in a graph. In other words, we intend to concentrate on the exclusive variable of the Bartholdi zeta function. To reach this aim, we ignore the other variable which is common with the Ihara zeta function. Moreover, the other advantage of omitting this variable lies in the fact that the determinant expression of the reduced Bartholdi zeta function can be computed in polynomial time.

The rest of the paper is organized as follows. In the next section, we establish the definition of the Bartholdi zeta function and present some results and preliminaries which are necessary for the remainder of the paper. In Section 3, the *semi principle minors* of the reduced Bartholdi zeta function are defined and their properties are studied. Section 4 presents the relation between the structure of graphs and the first five coefficients of the reduced Bartholdi zeta function. Finally, in Section 5, we achieve a general formula for calculating all the coefficients of the reduced Bartholdi zeta function.

## 2. Preliminaries

In this section, we state some definitions and theorems which are useful in the sequel of the paper.

Formally, a digraph $D$ is an ordered pair $(V(D), A(D))$ consisting of a non-empty set $V(D)$ of vertices, and a set $A(D) \subseteq V(D) \times V(D)$ of arcs. For any arc $a = (u, v)$, let $o(a) = u$ be the origin of $a$, and $t(a) = v$ be its terminus. Also the inverse of the arc $a = (u, v)$ obtained by switching the origin and terminus of $a$, is denoted by $a^{-1} = (v, u)$. A path $P$ of length $t$ in $D$ is a finite non-empty sequence $(v_0, a_1, v_1, \ldots, a_t, v_t)$, whose terms are alternately vertices and arcs, such that $a_i = (v_{i-1}, v_i)$, for $1 \leqslant i \leqslant t$. For abbreviation, we write $P = (a_1, \ldots, a_t)$.

We say that $P$ has a *backtracking* or *bump* if $a_{i+1} = a_i^{-1}$, for some $1 \leqslant i \leqslant t-1$. Moreover, if $v_0 = v_t$ then $P$ is called a cycle. Also, a cycle with no backtracking is called a *backtrack-less* cycle. Let $B^r$ be the cycle obtained by going $r$ times around a cycle $B$, which is called $r-$multiple of cycle $B$. The cycle $C$ is called *primitive* if it is not constructed by $r-$multiple of some other cycle $B$, i.e. for a natural number $r \geqslant 2$, $C \neq B^r$. Next, two cycles $C_1 = (a_1, a_2, \ldots, a_m)$ and $C_2 = (b_1, b_2, \ldots, b_m)$ are called *equivalent* if $b_j = a_{j+l}$, for a fixed number $l$ and for all $j$, where the subscripts are considered modulo $m$. Hence, the *equivalence class* of $C_1$ is

$$[C_1] = \{(a_1, a_2, \ldots, a_m), (a_2, a_3, \ldots, a_m, a_1), \ldots, (a_m, a_1, \ldots, a_{m-1})\}.$$

Now, consider a simple undirected and finite graph $G = (V(G), E(G))$ having $n$ vertices and $m$ edges, where $V(G)$ and $E(G)$ denote the vertex set and the edge set of $G$, respectively. Associated with the graph $G$, one may correspond a *symmetric digraph*, denoted by $\mathcal{D}(G)$ with $V(\mathcal{D}(G)) = V(G)$ and $\mathcal{E}(G) = A(\mathcal{D}(G)) = \{(u,v),(v,u)|uv \in E(G)\}$.

Recall that the oriented line graph of a digraph $D$, denoted by $\mathcal{L}(\mathcal{D}(G))$, is a digraph with vertex set $V(\mathcal{L}(\mathcal{D}(G))) = \mathcal{E}(G)$ and $(a_i, a_j) \in A(\mathcal{L}(\mathcal{D}(G)))$, if in the graph $\mathcal{D}(G)$, $t(a_i) = o(a_j)$ and $a_j \neq a_i^{-1}$, where $a_i, a_j \in \mathcal{E}(G)$.

The Ihara zeta function of an undirected graph $G$ is a function of a complex variable $t$, given by

$$Z_G(t) = \prod_{[C]}(1 - t^{l(C)})^{-1}, \tag{1}$$

where $[C]$ is the equivalence class of primitive backtrack-less cycles and the $l(C)$ is the length of the cycle $C$. For the matter of convergence of $Z_G(t)$, $|t|$ must be sufficiently small (see [6]).

Hashimoto generalized the determinant expression for the Ihara zeta function of an irregular graph by using the adjacency matrix of the oriented line graph (see [4]).

**Theorem 2.1.** *(Hashimoto). If $\boldsymbol{T}$ is the adjacency matrix of the oriented line graph $\mathcal{L}(\mathcal{D}(G))$, then*

$$Z_G(t)^{-1} = det(\boldsymbol{I}_{2m} - t\boldsymbol{T}), \tag{2}$$

*where $\boldsymbol{T} = (\boldsymbol{T}_{ij})$ can be defined as*

$$\boldsymbol{T}_{ij} = \begin{cases} 1 & \text{if } t(a_i) = o(a_j) \text{ and } a_j \neq a_i^{-1}, \\ 0 & \text{otherwise.} \end{cases}$$

In the digraph $\mathcal{D}(G)$, we say that a backtrack path $P = (a_1, \ldots, a_t)$ has a *bump* at $t(a_i)$ if $a_{i+1} = a_i^{-1}$ for some $1 \leqslant i \leqslant t-1$. The *cyclic bump count*, $cbc(C)$, is defined as the number of bumps in the cycle $C = (a_1, a_2, \ldots, a_t)$ and is given by

$$cbc(C) = |\{i = 1, \ldots, t | a_{i+1} = a_i^{-1}\}|. \tag{3}$$

Then the Bartholdi zeta function of $G$ is defined to be a function of $u, t \in \mathbb{C}$, with $|u|, |t|$ sufficiently small, given by

$$\zeta_G(u,t) = \zeta(u,t) = \prod_{[C]}(1 - u^{cbc(C)}t^{l(C)})^{-1}, \tag{4}$$

where $[C]$ runs over all equivalence classes of primitive cycles $C$ of $G$ (see [1]).

In the same paper, the determinant expression of the Bartholdi zeta function of $G$ is given as the following theorem.

**Theorem 2.2.** *(Bartholdi). Let $G$ be a connected graph with $n$ vertices and $m$ unoriented edges. Then the reciprocal of the Bartholdi zeta function of $G$ is given by*

$$\zeta(u,t)^{-1} = (1-(1-u)^2 t^2)^{m-n} det(\boldsymbol{I} - t\boldsymbol{A} + (1-u)(\boldsymbol{D} - (1-u)\boldsymbol{I})t^2), \tag{5}$$

*where $\boldsymbol{A}$ is the adjacency matrix of the graph $G$, and $\boldsymbol{D} = (d_{ij})$ is the degree matrix with $d_{ii} = deg(v_i)$, $v_i \in V(G)$ and for all $i \neq j$; $d_{ij} = 0$.*

For any $a_i, a_j \in \mathcal{E}(G)$, two $2|E| \times 2|E|$ matrices $\boldsymbol{B} = (\boldsymbol{B}_{ij})$, and $\boldsymbol{J} = (\boldsymbol{J}_{ij}) = \begin{bmatrix} 0 & \boldsymbol{I}_{|E|} \\ \hline \boldsymbol{I}_{|E|} & 0 \end{bmatrix}$ are defined as follows:

(6) $$\boldsymbol{B}_{ij} = \begin{cases} 1 & \text{if } t(a_i) = o(a_j), \\ 0 & \text{otherwise,} \end{cases} \quad \boldsymbol{J}_{ij} = \begin{cases} 1 & \text{if } a_j = a_i^{-1}, \\ 0 & \text{otherwise.} \end{cases}$$

Based on the definitions of $\boldsymbol{B}$ and $\boldsymbol{J}$, the determinant of the Bartholdi zeta function can also be written as (see [1])

(7) $$\zeta(u,t)^{-1} = \det(\boldsymbol{I}_{2m} - (\boldsymbol{B} - (1-u)\boldsymbol{J})t).$$

Furthermore, Horton in [5] stated some properties of the structure of the matrix $\boldsymbol{T}$, which come as the following remark:

**Remark 2.3.** *(Horton) Consider $\boldsymbol{T}$ as a block matrix $\begin{bmatrix} A & B \\ \hline C & D \end{bmatrix}$ where $A, B, C, D$ are $|E| \times |E|$ matrices with the following properties:*
*(i) $B = B^T, C = C^T$, and $D = A^T$.*
*(ii) The diagonal entries of $B$ and $C$ are zeros.*
*(iii) The diagonal entries of $A$ and $D$ are zeros whenever the graph contains no loops.*
*(iv) $\boldsymbol{T}^T = \boldsymbol{J}\boldsymbol{T}\boldsymbol{J}$.*
*(v) $A + B + C + A^T$ is the adjacency matrix of the line graph of $G$.*

The line graph of $G$ denoted by $L(G)$, is constructed by replacing each edge of $G$ with a vertex, where vertices of $L(G)$ are adjacent if and only if their corresponding edges in $G$ are incident to a same vertex.

Properties (i) and (ii) are explained by Stark and Terras [14]. As expressed in [5], property (iv) follows from property (i). Also properties (iii) and (v) are obtained from the definition of $\boldsymbol{T}$.

**Remark 2.4.** *Let $G$ be an undirected graph with the edge set $E(G) = \{e_1, \ldots, e_m\}$. Arbitrarily orient its edges and label all $2m$ arcs of its symmetric digraph $\mathcal{D}(G)$ as $\mathcal{E}(G) = \{a_1, \ldots, a_m, a_{m+1} = a_1^{-1}, \ldots, a_{2m} = a_m^{-1}\}$. Furthermore, let $a_{i'} = a_i^{-1}$ and $a_i = a_{i'}^{-1}$ for all $i$, $1 \leqslant i \leqslant m$. Hence, the labels of the arcs are changed to $\mathcal{E}(G) = \{a_1, \ldots, a_m, a_{1'}, \ldots, a_{m'}\}$. Since $\boldsymbol{T} = (c_{ij})$ is a block matrix, the labeling of its blocks based on the above, is obtained as $A = (c_{ij})_{1 \leqslant i \leqslant m; 1 \leqslant j \leqslant m}$, $B = (c_{ij})_{1 \leqslant i \leqslant m; 1' \leqslant j \leqslant m'}$, $C = (c_{ij})_{1' \leqslant i \leqslant m'; 1 \leqslant j \leqslant m}$ and $D = (c_{ij})_{1' \leqslant i \leqslant m'; 1' \leqslant j \leqslant m'}$.*

Before stating our results in the next sections, we need to derive the modified version of the Bartholdi zeta function. Therefore, we provide the following lemma and then define a new Bartholdi zeta function called the reduced Bartholdi zeta function.

**Lemma 2.5.** *Let $\boldsymbol{T}$ be the adjacency matrix of the oriented line graph of $G$. Then the determinant expression of the reciprocal of the Bartholdi zeta function can be given by*

(8) $$\zeta(u, \frac{1}{x})^{-1} = \frac{1}{x^{2m}} \det(x\boldsymbol{I}_{2m} - (\boldsymbol{T} + u\boldsymbol{J})).$$

*Proof.* Using the definitions of $\boldsymbol{B}$ and $\boldsymbol{J}$, clearly $\boldsymbol{T} = \boldsymbol{B} - \boldsymbol{J}$. Thus, equation (7) can be written as

(9) $$\zeta(u,t)^{-1} = \det(\boldsymbol{I}_{2m} - (\boldsymbol{T} + u\boldsymbol{J})t).$$

Then, by substituting $t = \frac{1}{x}$, the result follows. □

By factoring $\frac{1}{x^{2m}}$ from equation (8), we can write the reciprocal Bartholdi zeta function $\zeta(u, \frac{1}{x})^{-1}$ as $\chi(u, v) = \det(v\boldsymbol{I}_{2m} - (\boldsymbol{T} + u\boldsymbol{J}))$. If we set $v = 0$ in the polynomial $\chi(u, v = 0)$, we can define the reduced Bartholdi zeta function as bellow.

**Definition 2.6.** *The polynomial of $\chi(u, v = 0)$ is called the **reduced Bartholdi zeta function** and is given by*

(10) $\quad \chi(u, v = 0) = (-1)^{2m} \det(\boldsymbol{T} + u\boldsymbol{J}) = d_0 u^{2m} + d_1 u^{2m-1} + d_2 u^{2m-2} + \cdots + d_{2m}.$

By a direct investigation, one can verify that all the coefficients of the reduced Bartholdi zeta function enumerate the backtrack paths of the graph and are quite different from the Ihara zeta function which confirms on the enumeration of backtrack less paths. In the next section, we demonstrate how the backtrack paths of the graph can be used in the construction of the coefficients of the reduced Bartholdi zeta function. Moreover, in the remaining parts of the paper, we try to show the relation between these coefficients called *Backtrack coefficients of the reduced Bartholdi zeta function* and the structure of graph $G$.

## 3. The properties of Semi principle minors of the reduced Bartholdi zeta function

According to Remark 2.4, let $S[\alpha; \alpha'] = S_{q \times q}$ denote the $q \times q$ submatrix of $\boldsymbol{T}$ with rows indexed by a subset $\alpha = \{i_1, \ldots, i_q\} \subseteq \{1, \ldots, m, 1', \ldots, m'\}$ (for $q \geq 2$) and columns indexed by corresponding inverse subset $\alpha' = \{i'_1, \ldots, i'_q\}$ such that $i_1 < \cdots < i_q$, for all $i, 1 \leq i \leq m$.

**Definition 3.1.** *The determinant of the above submatrix $S[\alpha; \alpha']$ of $\boldsymbol{T}$ is called the **semi principal minor** of $\boldsymbol{T}$, obtained by selecting a subset of rows and the same inverse subset of the columns as the above. Also the matrix $S$ will be called as the semi principal minor matrix.*

**Example 3.2.** *Let $\boldsymbol{T} = (c_{ij})$ by Remark 2.4. An example of semi principal minor of $\boldsymbol{T}$, $S'_{4 \times 4}$, can be shown as follows:*

$$S'[1, 3, 2', 4'; 1', 3', 2, 4] = \begin{bmatrix} 0 & c_{13'} & c_{12} & c_{14} \\ c_{31'} & 0 & c_{32} & c_{34} \\ c_{2'1'} & c_{2'3'} & 0 & c_{2'4} \\ c_{4'1'} & c_{4'3'} & c_{4'2} & 0 \end{bmatrix}$$

**Theorem 3.3.** *Let $\boldsymbol{T}$ and $\boldsymbol{J}$ be two $2m \times 2m$ square matrices associated with $\mathcal{L}(\mathcal{D}(G))$ and define $\boldsymbol{M} = \boldsymbol{T} + u\boldsymbol{J}$. Then, the reduced Bartholdi zeta function of $G$ is given by $\det(\boldsymbol{M}) = d_0 u^{2m} + d_1 u^{2m-1} + d_2 u^{2m-2} + \cdots + d_{2m}$.*

*Let $d_k$ be the coefficient of the variable $u^{2m-k}, 0 \leq k \leq 2m$. Then, the value of $(-1)^m d_k$ is obtained by the sum of those semi principal minors of $\boldsymbol{T}$ which have $k$ rows and $k$ columns.*

*Proof.* Let $\boldsymbol{M} = \left[ \begin{array}{c|c} A' & B' \\ \hline C' & D' \end{array} \right]$. Based on labeling of the matrix $\boldsymbol{T}$, the entries of blocks of the matrix $\boldsymbol{M}$ can be shown as follows:

$$(11) \quad \boldsymbol{M} = \left[\begin{array}{cccc|cccc} 0 & c_{12} & \cdots & c_{1m} & \mathbf{u} & c_{12'} & \cdots & c_{1m'} \\ c_{21} & 0 & \cdots & c_{2m} & c_{21'} & \mathbf{u} & \cdots & c_{2m'} \\ \vdots & & \ddots & \vdots & \vdots & & \ddots & \vdots \\ c_{m1} & c_{m2} & \cdots & 0 & c_{m1'} & c_{m2'} & \cdots & \mathbf{u} \\ \hline \mathbf{u} & c_{1'2} & \cdots & c_{1'm} & 0 & c_{1'2'} & \cdots & c_{1'm'} \\ c_{2'1} & \mathbf{u} & \cdots & c_{2'm} & c_{2'1'} & 0 & \cdots & c_{2'm'} \\ \vdots & & \ddots & \vdots & \vdots & & \ddots & \vdots \\ c_{m'1} & c_{m'2} & \cdots & \mathbf{u} & c_{m'1'} & c_{m'2'} & \cdots & 0 \end{array}\right]$$

Then, by the property of the determinant, we have

$$\det(M) = (-1)^m \det\left[\begin{array}{c|c} B' & A' \\ \hline D' & C' \end{array}\right]$$

$$(12) \quad = (-1)^m \det \left[\begin{array}{cccc|cccc} \mathbf{u} & c_{12'} & \cdots & c_{1m'} & 0 & c_{12} & \cdots & c_{1m} \\ c_{21'} & \mathbf{u} & \cdots & c_{2m'} & c_{21} & 0 & \cdots & c_{2m} \\ \vdots & & \ddots & \vdots & \vdots & & \ddots & \vdots \\ c_{m1'} & c_{m2'} & \cdots & \mathbf{u} & c_{m1} & c_{m2} & \cdots & 0 \\ \hline 0 & c_{1'2'} & \cdots & c_{1'm'} & \mathbf{u} & c_{1'2} & \cdots & c_{1'm} \\ c_{2'1'} & 0 & \cdots & c_{2'm'} & c_{2'1} & \mathbf{u} & \cdots & c_{2'm} \\ \vdots & & \ddots & \vdots & \vdots & & \ddots & \vdots \\ c_{m'1'} & c_{m'2'} & \cdots & 0 & c_{m'1} & c_{m'2} & \cdots & \mathbf{u} \end{array}\right]$$

Now, let

$$\left[\begin{array}{c|c} A' & B' \\ \hline C' & D' \end{array}\right] = (a_{ij})_{1 \leqslant i,j \leqslant 2m}.$$

Then, by the fundamental property of the determinant, we have

$$(13) \quad d_k = (-1)^m \sum_{i_1 < i_2 < \cdots < i_k} \det\left[\begin{array}{cccc} a_{i_1 i_1} & a_{i_1 i_2} & \cdots & a_{i_1 i_k} \\ a_{i_2 i_1} & a_{i_2 i_2} & \cdots & a_{i_2 i_k} \\ \vdots & & \ddots & \vdots \\ a_{i_k i_1} & a_{i_k i_2} & \cdots & a_{i_k i_k} \end{array}\right]$$

$$= (-1)^m \times (\textit{the sum of semi principal minors of } \boldsymbol{T}).$$

□

**Remark 3.4.** *Let $\beta$ be the selected indices of rows in $\boldsymbol{T}$ and $\beta'$ be the inverse indices of its columns, associated to a semi principal minor $S = (S_{ij})$.*

*Then $\forall (i,j,k) \in \beta$ and the corresponding $(i',j',k') \in \beta'$, we have*

$$(14) \quad \textit{if } S_{ij'} = 1, S_{kj'} = 1 \implies S_{ik'} = 1, S_{ki'} = 1.$$

*Afterwards, we call this relation as the **relative transitive rule**.*

Before calculating the coefficients of the reduced Bartholdi zeta function of $G$, we have to consider some important properties of *semi principal minors*.

Let $S$ be a $q \times q$ *semi principle minor matrix* of $\boldsymbol{T}$. Each matrix $S$ represents an adjacency matrix of the oriented line graph of a subgraph of $\mathcal{D}(G)$ by Definition 3.1. Since $T$ is the adjacency matrix of the oriented line graph $\mathcal{L}(\mathcal{D}(G))$, the inverse operations of constructing $\mathcal{L}(\mathcal{D}(G))$ from $T$ can be applied to $S$, which results the desired subgraph. Furthermore, suppose that the related subgraph of $S$ in $\mathcal{D}(G)$ has an isolated edge which is not connected to the other edges in the subgraph. Thus, in the $\mathcal{L}(\mathcal{D}(G))$, this isolated edge can be considered as an isolated vertex which leads to $\det(S) = 0$. By contradiction, we can see that if $S$ has non-zero determinant, then the related subgraph does not consist of any isolated edge.

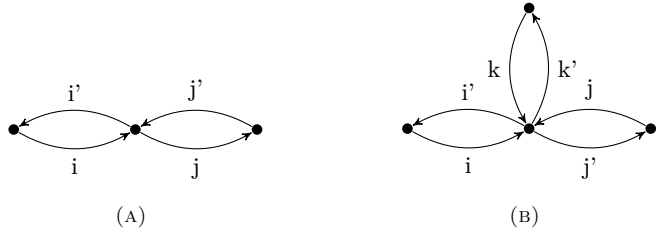

FIGURE 1. Two types of neighborhoods edges in a digraph.

According to Remark 2.4 and the definitions of $B$ and $J$, for any two adjacent edges, $a_i, a_j \in \mathcal{E}(G)$, we have $t(a_i) = o(a_j)$ and also $t(a_{j'}) = o(a_{i'})$ which are shown in Figure 1 (A). Thus, $S_{ij} = S_{j'i'}$, which shows that $S$ is a symmetric matrix. Also, on account of Remark 2.4 and Definition 3.1, the main diagonal of $S$ contains the entries $c_{ii'} \in T$ which are zeros, so the main diagonal of $S$ is always zero.

Moreover, by applying a similar argument, we can assert that if two edges $i$ and $k$ end to the origin of a common edge such as $j'$, then $t(a_i) = t(a_k) = o(a_{j'})$. Therefore, the symmetry of the $\mathcal{D}(G)$ will force them to have a relation between themselves as shown in Figure 1 (B). That is, for any matrix $S$, the *relative transitive rule* as "if $S_{ij'} = 1, S_{kj'} = 1$ then $S_{ik'} = 1, S_{ki'} = 1$" holds.

It is easily seen that for a symmetric permutation matrix, there is an even number of non zero off diagonal elements. For odd $q$, since the main diagonal of $S$ is zero, it is impossible that $S$ can be a permutation matrix.

From the above arguments in connection with the *semi principal minor matrices*, now we can state their properties as the following theorem.

**Theorem 3.5.** *Let $S$ be a $q \times q$ semi principle minor matrix of $T$. Then the following properties for the matrix $S$ hold:*
  (1) *Each $S$ represents an adjacency matrix of the oriented line graph of a subgraph of $\mathcal{D}(G)$.*
  (2) *If $S$ has non-zero determinant, then the related subgraph has no single isolated edge.*
  (3) *$S$ is a symmetric matrix, and its main diagonal is zero.*
  (4) *For all matrices $S$, the relative transitive rule holds.*
  (5) *If $q$ is odd, then $S$ never constructs a permutation matrix.*

4. REPRESENTATION OF THE COEFFICIENTS OF THE REDUCED BARTHOLDI ZETA FUNCTION

Let $G$ be a simple graph and $\chi(u, v = 0)$ be its reduced Bartholdi zeta function as equation (10). Then, for calculating the coefficients $d_0, d_1, d_2, d_3$ and $d_4$ of the reduced Bartholdi zeta function, we obtain the following results.

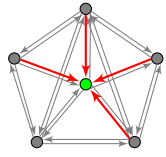

FIGURE 2. The sink star subgraph with 4 leaves.

**Lemma 4.1.** *The first backtrack coefficient of the reduced Bartholdi zeta function, $d_0$, has the value of $(-1)^m$ and the second coefficient, $d_1$, is equal to 0.*

*Proof.* By applying Theorem 3.3 to the matrix $\boldsymbol{M} = \boldsymbol{T} + u\boldsymbol{J}$, for calculating the coefficient $d_0$, all $2m$ variables $u$ must appear in *semi principle minors*. Therefore, the coefficient of $u^{2m}$ must be equal to $(-1)^m$ by which the result follows.

The same proof stands for $d_1 = 0$. While $2m-1$ variables of $u$ have been chosen, selecting the $1 \times 1$ permutation matrix as a *semi principle minor matrix* is always impossible. Thus, we have $(-1)^m u^{2m-1} \times 0 = 0$ as claimed. □

**Lemma 4.2.** *The coefficient $d_2$ of the reduced Bartholdi zeta function is determined by $(-1)^{m+1} \times ($ the number of edges of line graph of $G$, $L(G))$.*

*Proof.* By Theorem 3.3, we can see that the coefficient $d_2$ is $(-1)^m \times$ (the sum of all *semi principal minors* of $\boldsymbol{T}$ with 2 rows and 2 columns). We consider all such possible minors. Up to reordering the edges, the only non-trivial *semi principal minor* with 2 rows and 2 columns that may arise is:

$$\left| \begin{array}{cc} 0 & 1 \\ 1 & 0 \end{array} \right|$$

However, each entry of value "1" is selected only once in any *semi principle minor*. Thus, the number of different *semi principle minor* that only differ in their indices, is half of the number of "1"s in $\boldsymbol{T}$, which by Remark 2.3 is equal to twice the number of edges in $L(G)$. Therefore, the number of distinct *semi principle minors* of the above form in $\boldsymbol{T}$ has to be equal to the number of edges in $L(G)$. Since all these *semi principle minors* of size $2 \times 2$ have the determinant value of $-1$, the result follows. □

**Definition 4.3.** *A **sink star subgraph** in digraph $\mathcal{D}(G)$ is said to be a directed star subgraph containing a sink node $v \in V$ such that $deg^-(v) \geqslant 2$, and $q$ leaves with $q \leqslant deg^-(v)$. A set of all possible sink star subgraphs with $q$ leaves in digraph $\mathcal{D}(G)$ is denoted by $\mathscr{S}_q$.*

In Figure 2, a sink star subgraph selected from the digraph $\mathcal{D}(G)$, with 4 leaves is shown in red.

**Lemma 4.4.** *The number of sink star subgraphs in a digraph $\mathcal{D}(G)$, $|\mathscr{S}_q|$, with $q$ leaves can be written as*

(15) $$|\mathscr{S}_q| = \sum_{\substack{i=1 \\ d(v_i) \geqslant q}}^{n} \binom{d(v_i)}{q},$$

*where $d(v_i) = deg(v_i)$ is the degree of $i$-th vertex of $G$.*

*Proof.* It is obvious from Definition 4.3. □

**Lemma 4.5.** *Let $S$ be a $q \times q$ matrix of semi principle minor extracted from the matrix $\boldsymbol{T}$ of $G$ which consists of at least one zero off-diagonal entry. Then the indices of its rows represent the label of edges of unconnected sink star subgraphs in digraph $\mathcal{D}(G)$. Moreover, the sum of all indegrees of their sink nodes is equal to $q$, that is*

$$q = \sum_{v \in \text{ sink nodes of } S} deg^-(v)$$

*In other case, if $S$ consists of all-ones entries except its main diagonal, then the indices of its rows represent the labels of edges of a single sink star subgraph in digraph $\mathcal{D}(G)$, denoted by $|\mathscr{S}_q|$.*

*Proof.* Let us consider the $i$-th row of $S$ which corresponds to an edge in digraph $\mathcal{D}(G)$. This row must have at least one nonzero element, conforming to the properties (2) of Theorem 3.5 which emphasizes that the subgraph associated to $S$ has no single isolated edge.

Suppose that two columns $j', k'$ in this row have the value of 1 and the others are zero. Based on Definition 3.1, the connections between these three edges, $i, j'$ and $k'$ become similar to Figure 1 (B). Also, in row $j$, any columns except $i'$ and $k'$ have the value of zero, and similar arguments can apply to row $k$.

By Remark 3.4 and the property (4) of Theorem 3.5, for every matrix $S$, the **relative transitive rule** is known to hold. Hence, this rule makes the indices of any rows which participate in a sink star subgraph, to be disjoint from the other sink star subgraphs associated to $S$.

As a result, all indices of rows in $S$ construct a set of unconnected sink star subgraphs in digraph $\mathcal{D}(G)$. This fact that each index can contribute exactly to one sink star subgraph shows that the sum of all indegrees of their sink nodes will be equal to the number of indices of rows in $S$, which takes the value of $q$.

For the case that $S$ consists of all-ones entries except its main diagonal, by Theorem 3.5, it can be concluded that the corresponding sink star subgraph of the matrix $S$ has $q$ leaves in which all of its edges are incident to a sink node of a single sink star subgraph, denoted by $\mathscr{S}_q$ in Definition 4.3. □

**Lemma 4.6.** *The coefficient $d_3$ of the reduced Bartholdi zeta function is given by $(-1)^m \times 2 \sum_{\substack{i=1 \\ d(v_i) \geqslant 3}}^{n} \binom{d(v_i)}{3}$, where $d(v_i)$ is the degree of $i$-th vertex of $G$ and $m$ is the number of edges in $G$.*

*Proof.* By Theorem 3.3, the coefficient $d_3$ is equal to $(-1)^m \times$ (the sum of the *semi principal minors* of $\boldsymbol{T}$ with 3 rows and 3 columns). Based on their properties stated in Theorem 3.5, some of the non-trivial *semi principal minors*, with 3 rows and 3 columns are listed:

$$\begin{vmatrix} 0 & 1 & 1 \\ 1 & 0 & 0 \\ 1 & 0 & 0 \end{vmatrix}, \begin{vmatrix} 0 & 1 & 0 \\ 1 & 0 & 0 \\ 0 & 0 & 0 \end{vmatrix}, \ldots, \begin{vmatrix} 0 & 1 & 1 \\ 1 & 0 & 1 \\ 1 & 1 & 0 \end{vmatrix}$$

The determinant of all $3 \times 3$ *semi principal minors* that may arise are zero except the last one, which is the only non-zero determinant with the value of 2. From Lemma 4.5, when $S$ consists of all ones except its main diagonal, $|\mathscr{S}_3|$ can be calculated by Lemma 4.4, and the result will be obtained. □

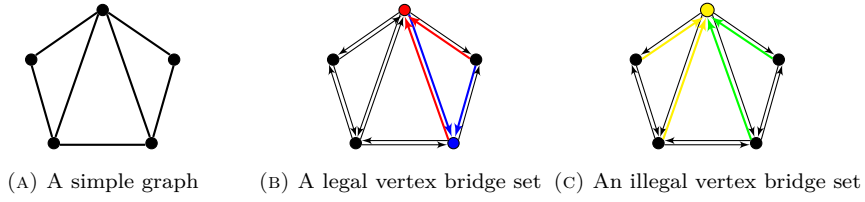

(A) A simple graph  (B) A legal vertex bridge set  (C) An illegal vertex bridge set

FIGURE 3. The types of vertex bridge set in the digraph $\mathcal{D}(G)$.

**Remark 4.7.** *Let a **vertex bridge** in $\mathcal{D}(G)$ be a sink star subgraph with 2 leaves. Each $v_i \in V$ ($1 \leqslant i \leqslant n$) with $d(v_i) \geqslant 2$ can be considered as the sink node of a vertex bridge. Then, all vertex bridges based on their sink nodes $v_i \in V$, can be divided into two sets, $\mathscr{B} = \{v_i \in V | d(v_i) = 2\}$ and $\mathscr{C} = \{v_i \in V | d(v_i) \geqslant 3\}$, based on the degrees of their sink nodes.*

*Let $C(v_i)$ be the number of vertex bridges that can be extracted from the vertex $v_i$ whose degree $d(v_i)$ in $G$ is at least 3. Then, by the sum over $C(v_i)$ for $d(v_i) \geqslant 3$, $|\mathscr{C}|$ is obtained and given by*

$$|\mathscr{C}| = \sum_{\substack{i=1 \\ d(v_i) \geqslant 3}}^{n} C(v_i) = \sum_{\substack{i=1 \\ d(v_i) \geqslant 3}}^{n} \binom{d(v_i)}{2}. \tag{16}$$

**Definition 4.8.** *The set of sink star subgraphs in the digraph $\mathcal{D}(G)$ with distinct sink nodes is called **legal set of sink star subgraphs**. Otherwise, the set of sink star subgraphs is called **illegal set of sink star subgraphs**, if it contains at least two sink star subgraphs in which their sink nodes are overlapped.*

According to Remark 4.7, we can define legal and illegal vertex bridge sets. In Figure 3, two examples of a legal set of sink star subgraphs and an illegal set of sink star subgraphs from the digraph $\mathcal{D}(G)$ with 5 vertices are shown.

**Lemma 4.9.** *Let $\mathcal{V}$ be the legal set of all pairs of vertex bridges in the digraph $\mathcal{D}(G)$. Then the size of $\mathcal{V}$ is given by*

$$|\mathcal{V}| = \binom{|\mathscr{C}| + |\mathscr{B}|}{2} - \sum_{\substack{i=1 \\ d(v_i) \geqslant 3}}^{n} \binom{C(v_i)}{2}. \tag{17}$$

*Proof.* The first term indicates the number of all pairs of vertex bridges in the union set $\mathscr{C} \cup \mathscr{B}$, which contains all possible vertex bridges, extracted from vertices in $G$, in which $d(v_i) = 2$ or $d(v_i) \geqslant 3$ ($1 \leqslant i \leqslant n$). The size of the illegal set of vertex bridges can be obtained by $C(v_i)$. Therefore, the illegal vertex bridges in which their sink nodes are overlapped, must be counted and written off from all possible pairs, which is equal to $\sum_{\substack{i=1 \\ d(v_i) \geqslant 3}}^{n} \binom{C(v_i)}{2}$. This is precisely the assertion of the lemma. □

**Theorem 4.10.** *The coefficient $d_4$ of the reduced Bartholdi zeta function is given as*

(18)
$$d_4 = (-1)^m \times \left( |\mathcal{V}| - 3 \sum_{\substack{i=1 \\ d(v_i) \geqslant 4}}^{n} \binom{d(v_i)}{4} \right)$$
$$= (-1)^m \times \left( \binom{|\mathscr{C}| + |\mathscr{B}|}{2} - \sum_{\substack{i=1 \\ d(v_i) \geqslant 3}}^{n} \binom{C(v_i)}{2} - 3 \sum_{\substack{i=1 \\ d(v_i) \geqslant 4}}^{n} \binom{d(v_i)}{4} \right).$$

*Proof.* All *semi principle minors* of size $4 \times 4$ of $T$, that satisfy the properties of Theorem 3.5, can be shown as follows:

$$\begin{vmatrix} 0 & 1 & 0 & 0 \\ 1 & 0 & 0 & 0 \\ 0 & 0 & 0 & 1 \\ 0 & 0 & 1 & 0 \end{vmatrix}, \begin{vmatrix} 0 & 0 & 1 & 0 \\ 0 & 0 & 0 & 1 \\ 1 & 0 & 0 & 0 \\ 0 & 1 & 0 & 0 \end{vmatrix}, \begin{vmatrix} 0 & 0 & 0 & 1 \\ 0 & 0 & 1 & 0 \\ 0 & 1 & 0 & 0 \\ 1 & 0 & 0 & 0 \end{vmatrix}, \begin{vmatrix} 0 & 1 & 1 & 1 \\ 1 & 0 & 1 & 1 \\ 1 & 1 & 0 & 1 \\ 1 & 1 & 1 & 0 \end{vmatrix}$$

It is obvious that the first three *semi principle minors* have the determinant of 1. By Theorem 3.5 (1) and Remark 4.7, it follows that these semi principle minors correspond to the legal sets containing pair of the vertex bridges in the digraph $\mathcal{D}(G)$ and for enumerating them we apply Lemma 4.9 to derive $|\mathcal{V}|$.

By Lemma 4.5, it can be concluded that the fourth *semi principle minor* corresponds to a class of sink star subgraphs with 4 leaves, $\mathscr{S}_4$ and has the determinant value of $-3$. Thus, $|\mathscr{S}_4|$, the last term in equation (18) is obtained by Lemma 4.4, which establishes the formula. $\square$

It is easily seen that if each pair of vertex bridges overlap at their sink nodes, they construct a sink star subgraph with 4 leaves which corresponds to the fourth *semi principle minor*. Note that the illegal sets of pair of vertex bridges are not considered as distinct pairs, and thus they are regarded as $\mathscr{S}_4$.

5. The structure of graphs and the reduced Bartholdi zeta function

Let $G$ be a simple undirected graph with $m$ edges and $\chi(u, v = 0)$ its reduced Bartholdi zeta function. For calculating the coefficients $d_k$ for all $k; 1 \leqslant k \leqslant 2m$ of the reduced Bartholdi zeta function, first we investigate the problem of a unique partition and demonstrate how it can be used in the purpose of achieving the main result.

**Definition 5.1.** *A **unique partition** of an integer $k$ is a non-increasing sequence of positive integers greater than 1, denoted by $c_1, c_2, \ldots, c_l$, such that $k = \sum_{i=1, c_i > 1}^{l} c_i$. Each $c_i, 1 \leqslant i \leqslant l$, is called a part of the partition. From now on, let $p(k)$ be the number of unique partitions of the integer $k$, and $\rho_k^i$ denote the i-th unique partition of integer $k$, for $1 \leqslant i \leqslant p(k)$.*

**Remark 5.2.** *Let $\rho_k^i$ be the i-th partition of integer $k$. Corresponding to each $c_i$, by Definition 4.3, a sink star graph $\mathscr{S}_{c_i}$ with $c_i$ leaves, $1 \leqslant i \leqslant l$, is associated. Furthermore, the collection of $\mathscr{S}_{c_1}, \ldots, \mathscr{S}_{c_i}, \ldots, \mathscr{S}_{c_l}$ in which $c_1 + \cdots + c_i + \cdots + c_l = k$, constructs a set of unconnected sink star subgraphs to which a semi principle minor of size $k$ can be assigned.*

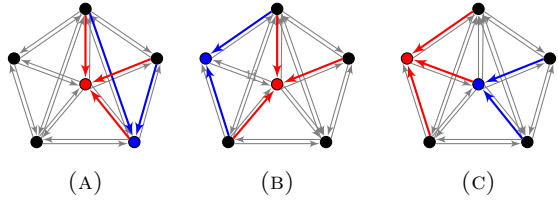

(A)     (B)     (C)

FIGURE 4. Some sink star subgraphs set contains $\{\mathscr{S}_2, \mathscr{S}_3\}$ in which $\mathscr{S}_2$ and $\mathscr{S}_3$ are shown in blue and red, respectively.

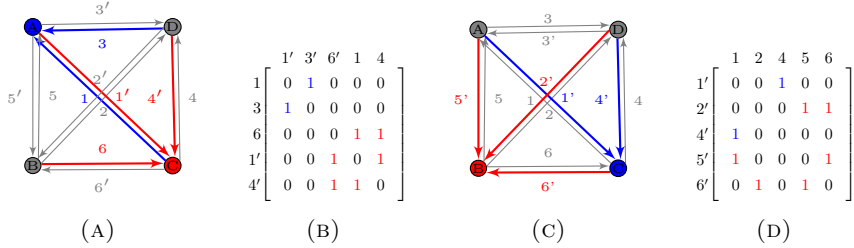

(A)     (B)     (C)     (D)

FIGURE 5. Some similar semi principle matrices extracted from isomorphic sink star subgraphs of an arbitrary digraph $\mathcal{D}(G)$.

By Definition 4.3, constructing $\mathscr{S}_{c_i}$ from the part $c_i$ is straightforward. The basic idea of linking between the partition problem and the *semi principle minor matrix* is that instead of looking at each $\mathscr{S}_{c_i}$ separately, we consider $\mathscr{S}_{c_1}, \ldots, \mathscr{S}_{c_i}, \ldots, \mathscr{S}_{c_l}$ as a set of unconnected sink star subgraphs of an arbitrary digraph. Using Lemma 4.5, the matrix of the corresponding *semi principle minor* can be obtained.

The relation between the unique partition problem of an integer $k$ and the corresponding unconnected sink star subgraphs $\mathscr{S}_{c_1}, \ldots, \mathscr{S}_{c_i}, \ldots \mathscr{S}_{c_l}$ with $c_1, \ldots, c_i, \ldots, c_l$ leaves, respectively, are shown in equation (19). Also, by Lemma 4.5, a *semi principle minor matrix* $S_{k \times k}$ can be assigned to them.

(19)
$$\begin{array}{c}
c_1 + \cdots + c_i + \cdots + c_l = k \quad \text{( Definition 5.1)} \\
\text{Remark 5.2} \downarrow \quad\quad \downarrow \quad\quad \downarrow \\
\{\mathscr{S}_{c_1}, \ldots, \mathscr{S}_{c_i}, \ldots, \mathscr{S}_{c_l}\} \\
\text{Definition 4.3} \downarrow \quad\quad \downarrow \quad\quad \downarrow \\
\bigstar_{c_1} \cdots \bigstar_{c_i} \cdots \bigstar_{c_l} \xrightarrow{\text{Lemma 4.5}} [\ S\ ]_{k \times k}
\end{array}$$

Three sets of unconnected sink star subgraphs of size 5 which contain $\{\mathscr{S}_2, \mathscr{S}_3\}$ are shown in the Figure 4.

**Remark 5.3.** *Suppose $A$ and $B$ are two $q \times q$ matrices of arbitrary semi principle minors. We say $A$ and $B$ are similar if $B = P^{-1}AP$, for some $q \times q$ permutation matrix $P$. Then the corresponding minors of matrices $A$ and $B$ are called **similar***

| k | $\rho_k^i$ | Determinant of $\mathscr{P}_k^i$ | *Some similar semi principle matrices* | *Related set of sink star subgraphs, $\mathscr{U}_k^i$* |
|---|---|---|---|---|
| 3 | (3) | 2 | $\begin{bmatrix}0&1&1\\1&0&1\\1&1&0\end{bmatrix}$ | 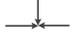 |
| 4 | (2,2) | 1 | $\begin{bmatrix}0&0&0&1\\0&0&1&0\\0&1&0&0\\1&0&0&0\end{bmatrix},\begin{bmatrix}0&0&1&0\\0&0&0&1\\1&0&0&0\\0&1&0&0\end{bmatrix},\begin{bmatrix}0&1&0&0\\1&0&0&0\\0&0&0&1\\0&0&1&0\end{bmatrix},\ldots$ | 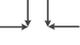 |
| 4 | (4) | -3 | $\begin{bmatrix}0&1&1&1\\1&0&1&1\\1&1&0&1\\1&1&1&0\end{bmatrix}$ | 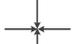 |
| 5 | (2,3) | -2 | $\begin{bmatrix}0&0&1&0&0\\0&0&0&1&1\\1&0&0&0&0\\0&1&0&0&1\\0&1&0&1&0\end{bmatrix},\begin{bmatrix}0&0&0&0&1\\0&0&1&1&0\\0&1&0&1&0\\0&1&1&0&0\\1&0&0&0&0\end{bmatrix},\begin{bmatrix}0&0&1&0&1\\0&0&0&1&0\\1&0&0&0&1\\0&1&0&0&0\\1&0&1&0&0\end{bmatrix},\begin{bmatrix}0&1&1&0&0\\1&0&1&0&0\\1&1&0&0&1\\0&0&0&0&1\\0&0&0&1&0\end{bmatrix},\begin{bmatrix}0&0&1&0&0\\0&0&0&1&1\\1&0&0&0&0\\0&1&0&0&1\\0&1&0&1&0\end{bmatrix},\ldots$ | 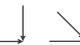 |
| 5 | (5) | 4 | $\begin{bmatrix}0&1&1&1&1\\1&0&1&1&1\\1&1&0&1&1\\1&1&1&0&1\\1&1&1&1&0\end{bmatrix}$ | 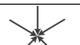 |

TABLE 1. The $i$-th partition of $k$, $\rho_k^i$, its prototype matrix, $\mathscr{P}_k^i$, and its related set of sink star subgraphs, $\mathscr{U}_k^i$, for $k = 3, \ldots, 5$.

***semi principle minors***. *Furthermore, the similar semi principle minors have the same determinant.*

**Lemma 5.4.** *Two sets of unconnected sink star subgraphs $H, F$ are isomorphic if and only if their semi principle minors are similar.*

*Proof.* Clearly, the result follows by the property of the isomorphic graphs and Remark 5.3. □

**Remark 5.5.** *Let $\rho_k^i$ be the $i$-th partition of integer $k$. Corresponding to each $c_i$, by Definition 4.3, a sink star graph $\mathscr{S}_{c_i}, 1 \leqslant i \leqslant l$, is associated. Let $\mathscr{U}_k^i$ be the set of sink star subgraphs, $\mathscr{S}_{c_1}, \ldots, \mathscr{S}_{c_i}, \ldots, \mathscr{S}_{c_l}$ corresponding to $\rho_k^i$ in the digraph $\mathcal{D}(G)$.*

*Each $\mathscr{U}_k^i$ can be associated to some distinct semi principle minors with a different form of size $k$ in which all pairs of them are similar. On the contrary, a set of similar semi principle minors of size $k$ represents a set of isomorphic $\mathscr{U}_k^i$.*

In Figure 5, two *similar semi principle minors* extracted from two isomorphic subgraphs are shown.

**Definition 5.6.** *Suppose $\rho_k^i$, for all $k, 1 \leqslant i \leqslant p(k)$. We choose one matrix of similar semi principle minors as a* **prototype matrix** *for $\rho_k^i$, called $\mathscr{P}_k^i$, which exhibits a set of isomorphic subgraphs.*

In Table 1, some prototype matrices $\mathscr{P}_k^i$ for the $i$-th partition of an integer $k$, $\rho_k^i$, for $k = 3, \ldots, 5$ are shown. Also, the value of determinant of the prototype matrices, their similar semi principle matrices and their related sink star subgraphs, $\mathscr{U}_k^i$, are presented.

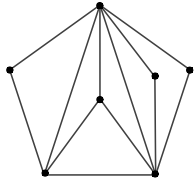

Figure 6. A simple graph $G$ of order 7.

**Theorem 5.7.** *Let $\mathscr{U}_k^i$ be the set of sink star subgraphs corresponding to $\rho_k^i$ in $\mathcal{D}(G)$. Then the k-th coefficient, $d_k$, of the term $u^{2m-k}$ in the reduced Bartholdi zeta function of the graph $G$ is given by*

$$d_k = (-1)^m \sum_{i=1}^{p(k)} |\mathscr{U}_k^i| \times \boldsymbol{det}(\mathscr{P}_k^i), \tag{20}$$

*where $\mathscr{P}_k^i$ is the prototype matrix of $\rho_k^i$ and $p(k)$ denotes the number of unique partitions of the integer $k$.*

*Proof.* By Theorem 3.3, for each $k \in \{0,1,2,...,2m\}$, the coefficient $(-1)^m d_k$ is the sum of all the *semi principal minors* of the matrix $\boldsymbol{T}$, with $k$ rows and $k$ columns. On the other hand, by Remark 5.2, each *semi principle minor* of size $k$ can be taken into account as a set of $\mathscr{S}_{c_1}, \ldots, \mathscr{S}_{c_i}, \ldots, \mathscr{S}_{c_l}$ where the sum of their indices is equal to $k$. Enumerating all these sets leads to the problem of partitioning.

From Lemma 4.5 and Remark 5.5, we can see that the sink star subgraphs of different $\rho_k^i$ will cover all the types of *semi principal minors* of size $k$. Moreover, using Definition 4.8, the determinant of the prototype matrix of $\rho_k^i$, $\mathscr{P}_k^i$, must be multiplied by the number of sink star subgraphs of $\rho_k^i$ to obtain the result as claimed. □

We give an example to show how we can apply Theorem 5.7 to a graph for the purpose of finding a coefficient of its reduced Bartholdi zeta function.

**Example 5.8.** *Let us consider the simple graph $G$ with 7 vertices and 12 edges as shown in Figure 6. Suppose that we want to obtain the coefficient $d_5$ of $u^{19}(u^{2m-5})$ as defined in the reduced Bartholdi zeta function of $G$, by the way proposed in this paper. The complete polynomial of the reduced Bartholdi zeta of this graph is given as follows:*

$$u^{24} - 37u^{22} + 70u^{21} + 435u^{20} - 1708u^{19} - 183u^{18} + 11646u^{17} - 21726u^{16}$$
$$- 6800u^{15} + 80478u^{14} - 111316u^{13} + 3094u^{12} + 171768u^{11} - 209950u^{10}$$
$$+ 49164u^9 + 134181u^8 - 163728u^7 + 66487u^6 + 23310u^5 - 44185u^4$$
$$+ 26068u^3 - 8475u^2 + 1526u - 120.$$

First, we have to find the parts of the unique partition of integer 5 such as $(2,3),(5)$. According to Definition 5.1, let $\rho_5^1 = (2,3), \rho_5^2 = (5)$.

The number of all possible sink star subgraphs of $\rho_5^1$ such as $\mathscr{S}_2, \mathscr{S}_3$ and $\rho_5^2$ such as $\mathscr{S}_5$ must be determined.

In order to find all legal pairs of $\mathscr{S}_2, \mathscr{S}_3$ in $\mathcal{D}(G)$, we calculate the value of $\binom{d(v_i)}{2}\binom{d(v_j)}{3}, 1 \leqslant i,j \leqslant 7$. The vertices $v_i, 1 \leqslant i \leqslant 7$, are selected as the sink nodes of $\mathscr{S}_2$ and the vertices $v_j, 1 \leqslant j \leqslant 7$, are selected as the sink nodes of $\mathscr{S}_3$. By considering the number of degrees and by Definition 4.8, we conclude that the pairs of $\binom{3}{2}\binom{3}{3}, \binom{4}{2}\binom{4}{3}, \binom{5}{2}\binom{5}{3}, \binom{6}{2}\binom{6}{3}$ are illegal sets, since the sink nodes of $\mathscr{S}_2, \mathscr{S}_3$ are overlapped and we set the number of them as zero. Now we have to sum over all the combinations of legal pairs listed in bellow in which contains 868 elements.

$$\begin{array}{cccc}
3\binom{2}{2}\binom{3}{3} & 3\binom{2}{2}\binom{4}{3} & 3\binom{2}{2}\binom{5}{3} & 3\binom{2}{2}\binom{6}{3} \\
0 & \binom{3}{2}\binom{4}{3} & \binom{3}{2}\binom{5}{3} & \binom{3}{2}\binom{6}{3} \\
\binom{4}{2}\binom{3}{3} & 0 & \binom{4}{2}\binom{5}{3} & \binom{4}{2}\binom{6}{3} \\
\binom{5}{2}\binom{3}{3} & \binom{5}{2}\binom{4}{3} & 0 & \binom{5}{2}\binom{6}{3} \\
\binom{6}{2}\binom{3}{3} & \binom{6}{2}\binom{4}{3} & \binom{6}{2}\binom{5}{3} & 0
\end{array}$$

Moreover, for calculating $\rho_5^2$, the number of sink star subgraphs with 5 leaves, $\mathscr{S}_5$, must be determined. By Lemma 4.4, we know that

$$|\mathscr{S}_5| = \sum_{\substack{i=1 \\ d(v_i) \geqslant 5}}^{7} \binom{d(v_i)}{5} = 7.$$

As a result, considering the prototype matrices $\mathscr{P}_5^1$ and $\mathscr{P}_5^2$ extracted from Table 1 and by applying Theorem 5.7, we have

$$d_5 = (-1)^{12}((868 \times \det(\mathscr{P}_5^1)) + (7 \times \det(\mathscr{P}_5^2)))$$
$$= 868 \times (-2) + 7 \times 4 = -1708,$$

which is the coefficient $d_5$ of $u^{19}$ in the reduced Bartholdi zeta function.

## References

[1] Bartholdi, L., 2000. 'Counting paths in graphs', *arXiv preprint math/0012161*.
[2] Bass, H., 1992. 'The Ihara-Selberg zeta function of a tree lattice', *International Journal of Mathematics* **3**(06), 717–797.
[3] Cooper, Y., 2009. 'Properties determined by the Ihara zeta function of a graph', *Electron. J. Combin* **16**(1).
[4] Hashimoto, K.-i., 1992. 'Artin type L-functions and the density theorem for prime cycles on finite graphs', *International Journal of Mathematics* **3**(06), 809–826.
[5] Horton, M. D., 2006. 'Ihara zeta functions of irregular graphs'.
[6] Ihara, Y., 1966. 'On discrete subgroups of the two by two projective linear group over $\mathfrak{p}$-adic fields', *Journal of the Mathematical Society of Japan* **18**(3), 219–235.
[7] Kotani, M. and Sunada, T., 2000. 'Zeta functions of finite graphs'.
[8] Mizuno, H. and Sato, I., 2007. 'A new Bartholdi zeta function of a digraph', *Linear algebra and its applications* **423**(2), 498–511.
[9] Mizuno, H. and Sato, I., 2014. 'Some weighted Bartholdi zeta function of a digraph', *Linear Algebra and its Applications* **445**, 1–17.
[10] Oren, I., Godel, A. and Smilansky, U., 2009. 'Trace formulae and spectral statistics for discrete Laplacians on regular graphs (i)', *Journal of Physics A: Mathematical and Theoretical* **42**(41), 415101.

[11] Sato, I., 2013. 'A generalized Bartholdi zeta function for a hypergraph', *Far East Journal of Mathematical Sciences* **78**(1), 93.
[12] Sato, I., Mitsuhashi, H. and Morita, H., 2015. 'A generalized Bartholdi zeta function for a general graph', *Linear and Multilinear Algebra* pp. 1–18.
[13] Scott, G. and Storm, C., 2008. 'The coefficients of the Ihara zeta function', *Involve, a Journal of Mathematics* **1**(2), 217–233.
[14] Stark, H. M. and Terras, A. A., 1996. 'Zeta functions of finite graphs and coverings', *Advances in Mathematics* **121**(1), 124–165.
[15] Storm, C., 2007*a*. Extending the Ihara-Selberg zeta function to hypergraphs, PhD thesis, DARTMOUTH COLLEGE Hanover, New Hampshire.
[16] Storm, C. K., 2007*b*. 'Some graph properties determined by edge zeta functions', *arXiv preprint arXiv:0708.1923* .
[17] Terras, A. A. and Stark, H. M., 2007. 'Zeta functions of finite graphs and coverings, part III', *Advances in Mathematics* **208**(1), 467–489.

DEPARTMENT OF MATHEMATIC AND COMPUTER SCIENCE, AMIRKABIR UNIVERCITY OF TEHRAN, IRAN

*E-mail address*: `mstahaei@aut.ac.ir, nhashemi@aut.ac.ir`